\documentclass[11pt,leqno]{article}
\usepackage{amssymb}
\usepackage{amscd}
\usepackage{amscd}
\usepackage{amsmath}
\usepackage{amsthm}

\newcommand{\Q}{\mathbb Q}
\newcommand{\Z}{\mathbb Z}
\newcommand{\C}{\mathbb C}
\newcommand{\R}{\mathbb R}
\newcommand{\A}{\mathbb A}

\newcommand{\hH}{\mathcal H}
\newcommand{\Oo}{\mathcal O}
\newcommand{\nN}{\mathcal N}

\makeatletter
\def\1@section{\@tocline{1}{4pt}{1pc}{}{}}
\def\1@subsection{\@tocline{2}{0pt}{2pc}{5pc}{}}
\makeatother

\begin{document}
\title{\bf MODULAR CURVES, MODULAR SURFACES, AND MODULAR FOURFOLDS}
\author{{\bf Dinakar Ramakrishnan}
\\{}\\{\bf To Jacob Murre}}
\date{}

\setcounter{page}{1} \maketitle

\section{Introduction}

\medskip

We begin with some general remarks. Let $X$ be a smooth projective
variety of dimension $n$ over a field k. For any positive integer $p
<n$, it is of interest to understand, modulo a natural equivalence,
the algebraic cycles $Y=\sum_j m_j Y_j$ lying on $X$, with each
$Y_j$ closed and irreducible of codimension $p$, together with
codimension $p+1$ algebraic cycles $Z_j = \sum_i r_{ij}Z_{ij}$ lying
on $Y_j$, for all $j$. There is a natural setting in which to study
such a chain $(X\supset Y_j \supset Z_{ij})_{ij}$ of cycles, namely
when the following hold:
\begin{enumerate}
\item[(a)] Each $Z_j$ is, as a divisor on $Y_j$, linearly equivalent
to zero, i.e.,
of the form div$(f_j)$ for a function $f_j$ on $Y_j$;
\item[(b)] The formal sum $\sum_j m_j Z_j$ is zero as a codimension
$p+1$ cycle on $X$.
\end{enumerate}

\medskip

Those satisfying (a), (b) form a group ${\mathcal Z}^{p+1}(X,1)$. An
easy way to construct elements of this group is to take a
codimension $p-1$ subvariety $W$ of $X$, with a pair of (non-zero)
functions $(\varphi, \psi)$ on $W$, and take the formal sum $\sum_j
(Y_j, T_{Y_j}(\varphi,\psi))$, where $\{Y_j\}$ is the finite set of
codimension $p$ subvarieties where $\varphi$ or $\psi$ has a zero or
a pole, and
$$
T_{Y_j}(\varphi,\psi) \, = \, (-1)^{{\rm ord}_j(\varphi){\rm
ord}_j(\psi)}\left(\frac{\varphi^{{\rm ord}_j(\psi)}}{\psi^{{\rm
ord}_j(\varphi)}}\right)\vert_{Y_j},
$$
the {\it Tame symbol} of $(\varphi, \psi)$ at $Y_j$, where ord$_j$
denotes the order at $Y_j$. It is a fact that $\sum_j
T_{Y_j}(\varphi,\psi)$ is zero as a codimension $p+1$ cycle on $X$.
Let $CH^{p+1}(X,1)$ denote the quotient of ${\mathcal Z}^{p+1}(X,1)$
by the subgroup generated by such elements. This is a basic example
of Bloch's higher Chow group ([B]). For any abelian group $A$, let
$A_\Q$ denote $A \otimes \Q$. Then it may be worthwhile to point out
the isomorphisms
$$
CH^{p+1}(X,1)_\Q \, \simeq \, H_{\mathcal M}^{2p+1}(X,\Q(p+1)) \,
\simeq \, {\rm Gr}_\gamma^{p+1}K_1(X)\otimes \Q,
$$
where $H_{\mathcal M}^\ast(X,\Q(\ast\ast))$ denotes the bigraded
{\it motivic cohomology} of $X$, and Gr$_\gamma^r$ denotes the
$r$-th graded piece defined by the gamma filtration on $K_\ast(X)$.

There is another way to construct classes in this group, and that is
to make use of the product map
$$
CH^p(X) \otimes k^\ast \, \rightarrow \, CH^{p+1}(X,1),
$$
the image of this being generated by the {\it decomposable classes}
$\{(Y,\alpha)\}$, where each $Y$ is a codimension $p$ subvariety and
$\alpha$ a non-zero scalar.

\medskip

It should also be mentioned, for motivation, that when $k$ is a
number field, Beilinson predicts that the elements in
$CH^{p+1}(X,1)_\Q$ which come from a regular proper model $\mathcal
X$ over $\Z$ span a $\Q$-vector space of dimension equal to the
order of vanishing at $s=p$ of the associated $L$-function
$L(s,H^{2p}(X))$, which we will denote by $L(s,X)$ if there is no
confusion. By the expected functional equation, the order of pole at
$s=p+1$ of $L(s,X)$, denoted $r_{\rm an}(X)$, will be the difference
of the order of pole of the Gamma factor $L_\infty(s,X)$ at $s=p$
and the order of zero of $L(s,X)$ at $s=p$. The celebrated
conjecture of Tate asserts that $r_{\rm an}(X)$ is the dimension of
Im$\left(CH^p(X) \rightarrow H^{2p}_{\rm et}(X_{\overline \Q},
\Q_\ell)\right)$ for any prime $\ell$. One of the main objects here
is to sketch a proof of the Tate conjecture in codimension $2$ for
Hilbert modular fourfolds, and also deduce the Hodge conjecture
under a hypothesis.

\medskip

Going back to $CH^{p+1}(X,1)$, the first case of interest is when
$X$ is a surface and $p=1$. When $X$ is the Jacobian $J_0(37)$ of
the modular curve $X_0(37)$, Bloch constructed a non-trivial example
$\beta \in CH^2(X,1)_\Q$ by using the curve and the fact that
$J_0(37)$ is isogenous to a product of two elliptic curves over
$\Q$. This was generalized by Beilinson ([Be]; see also [Sch]) to a
product of two modular curves by going up to a (ramified) cover
$X_0(N) \times X_0(N)$ and by taking $\{Y_j\}$ to be the union of
the diagonal $\Delta$ and the curves $X_0(N) \times \{P\}$ and
$\{Q\} \times X_0(N)$, where $P, Q$ are cusps; the existence of the
functions $f_j$ came from the Manin-Drinfeld theorem saying that the
difference of any two cusps is torsion in the Jacobian. Later the
author generalized this ([Ra1,2]) to the case of Hilbert modular
surfaces $X$ by using a class of curves on $X$ called the
Hirzebruch-Zagier cycles, carefully chosen to have appropriate
intersection properties; in general these curves meet in CM points
or cusps.

\medskip

The second main goal of this article is to describe briefly the
ideas behind an ongoing project of the author involving the
construction of ($\Q$-rational) classes in $CH^3(X,1)_\Q$ for
certain modular fourfolds $X/\Q$. We will restrict ourselves to
Hilbert modular fourfolds defined by a biquadratic, totally real
field $F$. Very roughly, the basic idea is to use suitable
translates $Y_j$ of embedded Hilbert modular surfaces (coming from
the three quadratic subfields), choose $Z_j$ to be made up of
(translates of) embedded modular curves which are homologically
trivial, hence rationally trivial (as the $Y_j$ are simply
connected), and also use the fact that the Tate and Hodge
conjectures are known (by the author [Ra4]) for codimension $2$
cycles on $X$ (for square-free level), as well as the knowledge
({\it loc. cit.}) that a basis of $\Q$-rational cycles modulo
homological equivalence (in the middle dimension) on $X$, and on the
embedded Hilbert modular surfaces, is given by appropriate
Hirzebruch-Zagier cycles and their twists. We will first construct
decomposable classes and then indicate some candidates for the
indecomposable part.

\medskip

The ultimate goal is to understand this phenomenon for Shimura
varieties $X$, of which modular curves and Hilbert modular varieties
are examples. It is known that the most interesting (cuspidal) part
of the cohomology of $X$ is in the middle dimension, which leads us
to consider such $X$ of dimension $n=2m$ and take $p=m$. When there
are Shimura subvarieties $Y$ of $X$ of dimension $m$ with
$H^1(Y)=0$, like for Siegel modular varieties, one can hope to
construct promising classes in $CH^{m+1}(X,1)_\Q$. This will be
taken up elsewhere.

\medskip

This article is dedicated to Jaap Murre, from whom I have learnt a
lot over the years -- about algebraic cycles and about the
(conjectural) Chow-K\"unneth decompositions, though they exist for
simple reasons in the cases considered here. We have a long term
collaboration as well on the zero cycles on abelian surfaces. I
would also like to acknowledge a helpful conversation I had with
Spencer Bloch about $CH^\ast(X,1)$ in the modular setting (see
section 14). I thank the referee and Mladen Dimitrov for spotting
various typos on an earlier version, and for making suggestions for
improvement of exposition. Finally, I am pleased to acknowledge the
support of the National Science Foundation through the grant
DMS-0402044.

\bigskip

\section{Notations}

\medskip

Let $X$ be a smooth projective fourfold over a number field $k$.
Set:

\medskip

\begin{enumerate}
\item[]$V_B=H^4(X(\C),\Q)$ \item[]$Hg^2(X)= V_B\cap
H^{2,2}(X(\C))$ \item[]$r_{\rm Hg} = {\rm dim}_\Q Hg^2(X)$
\item[]${\mathcal G}_k = {\rm Gal}(\overline k/k)$
\item[]$\ell$: \, a prime
\item[]$V_\ell=H^4(X_{\overline \Q},\Q_\ell)$: ${\mathcal
G}_k$-module \item[]$r_{{\rm alg},k} = {\rm dim}\left({\rm
Im}(CH^2(X)_\Q\to V_\ell(2))\right)$
\item[]${\mathcal T}_{\ell,
k} = V_\ell(2)^{{\mathcal G}_k}$ \item[]$r_{\ell, k} = {\rm
dim}_{\Q_\ell} {\mathcal T}_{\ell, k}$
\item[]$S_{\rm bad}:= \, \{{\mathcal P} \, \vert
\, V_\ell \ne V_\ell^{I_{\mathcal P}}\} \cup\{{\mathcal P}\mid
\ell\}$
\item[]$S$: a finite
set of places $\supset S_{\rm bad}$
\item[]$Fr_{\mathcal P}$: \, geometric Frobenius at ${\mathcal P},
\, \forall {\mathcal P} \notin S_{\rm bad}$ \, with norm
$N({\mathcal P})$
\item[]$L(s,X) = \prod\limits_{\mathcal P \notin S} {\rm
det}(I-Fr_{\mathcal P}T \vert V_\ell)_{\vert_{T = N({\mathcal
P})^{-s}}}^{-1}$
\end{enumerate}

\medskip

By Deligne's proof of the Weil conjectures, the inverse roots of
$Fr_{\mathcal P}$ on $V_\ell$ are, for ${\mathcal P}\notin S$, of
absolute value $N({\mathcal P})^2$, implying that the $L$-function
$L(s,X)$ converges absolutely in $\{\Re(s)>3\}$. The boundary point
$s=3$, where $L(s,X)$ could be divergent, is called the Tate point.
The fourfolds of interest to us will admit meromorphic continuation
to the whole $s$-plane and satisfy a functional equation relating
$s$ to $5-s$. Put
$$
r_{{\rm an}, k} \, = \, -{\rm ord}_{s=3}L(s,X).
$$
Tate's conjecture is that this {\it analytic rank} equals the {\it
algebraic rank} $r_{{\rm alg},k}$ of codimension $2$ algebraic
cycles on $X$ modulo ($\ell$-adic) homological equivalence. It is
also expected that these two ranks are the same as the {\it
$\ell$-adic cycle rank} $r_{{\ell},k}$, and one always has $r_{{\rm
alg},k} \leq r_{{\ell},k}$.

\bigskip

\section{Hilbert Modular fourfolds}

\medskip

Let $K$ be a quartic, Galois, totally real number field with
embedding $K \hookrightarrow \R^4$ given by the archimedean places.
Fix a square-free ideal ${\nN}$ in the ring $\Oo_K$ of integers of
$K$, and write $\Gamma$ for the congruence subgroup
$\Gamma_1(\nN)\subset {\rm SL}(2, \Oo_K)$ of level $\nN$. Then there
is a natural embedding
$$
\Gamma \hookrightarrow {\rm SL}(2,\R)\times{\rm SL}(2,\R)\times{\rm
SL}(2,\R)\times{\rm SL}(2,\R), \, \gamma \to
(\gamma^\sigma)_{\sigma\in{\rm Hom}(K,\R)}.
$$
Using this one gets an action of $\Gamma$ on the four-fold product
of the upper half plane $\hH = {\rm SL}(2,\R)/{\rm SO}(2)$. The
quotient $Y=\Gamma\backslash{\hH}^4$ is a coarse moduli space of
polarized abelian fourfolds $A$ with $\Gamma$-structure, with
End$(A) \hookleftarrow \Oo \subset K$. It is a quasi-projective
variety, with Baily-Borel-Satake compactification $Y^\ast$, and a
smooth toroidal compactification $X: = \tilde Y =
Y\cup\tilde{Y}^\infty$, all defined over $\Q$.

For simplicity of exposition, we have used here the classical
formalism. Later, we will need to work with the adelic version
$S_{C_1(\nN)}$ relative to the standard compact open subgroup
$C_1(\nN)$ of $G(\A_{K,f})$, where $\A_{K,f}$ denotes the ring of
finite adeles of $K$; one has $C_1(\nN)\cap {\rm GL}(2,F) \, = \,
\Gamma_1(\nN)$. Moreover,
$$
S_{C_1(\nN)}(\C) \, = \, {\rm GL}(2, K)\backslash (\C-\R)^4 \times
{\rm GL}(2, \A_{K,f})/C_1(\nN),
$$
which is finitely connected, and $\Gamma_1(\nN)\backslash \hH^4$
occurs as an \'etale quotient of the connected component. (If the
object is to realize $\Gamma_1(\nN)\backslash \hH^4$ as {\it
exactly} the connected component, one needs to consider instead the
Shimura variety associated to the $\Q$-group $G$ with $G(\Q) =\{g\in
{\rm GL}(2,K)\, \vert \, {\rm det}(g)\in\Q^\ast\}$.) The Shimura
variety $S_{C_1(\nN)}$ is defined over $\Q$, and the same holds for
its Baily-Borel-Satake compactification $S_{C_1(\nN)}^\ast$. One can
also choose a smooth toroidal compactification
$X=\tilde{S}_{C_1(\nN)}$ over $\Q$. Again we will use $Y$, resp.
$\tilde Y^\infty$, to denote $S_{C_1(\nN)}$, resp. the boundary, so
that $X = Y\cup\tilde Y^\infty$.

\bigskip

\section{Results on Cycles of codimension $2$ on $X$}

\medskip

Let $X$ be the smooth toroidal compactification over $\Q$  of a
Hilbert modular fourfold of {\it square-free} level $\nN$ as above,
relative to a quartic Galois extension $K$ of $\Q$.

\medskip

\noindent{\bf Theorem A} ([Ra4]) \, \it \begin{enumerate}
\item[(i)]The Tate classes in $V_\ell(2)$ are algebraic. In fact,
$r_{{\ell},k}=r_{{\rm alg},k}=r_{{\rm an},k}$. \item[(ii)]If $\nN$
is a proper ideal, the Hodge classes in $V_B(2)$ are algebraic when
they are not pull-backs of classes from the full level, and
moreover, they are not all generated by intersections of divisors.
\end{enumerate}
\rm

\medskip

The next few sections indicate a proof of this, while at the same
time developing the theory and setting the stage for what is to come
afterwards.

\medskip

Now let $K$ be {\it biquadratic} so that Gal$(K/\Q) =
\{1,\sigma_1,\sigma_2, \, {\sigma_3=\sigma_1\sigma_2}\}$, with
$\sigma_j^2=1$ for each $j$. Let $F_j \subset K$ be the real
quadratic field obtained as the fixed field of $\sigma_j$. For every
$g\in $GL$_2^+(K)$, let $Y_{F_j,g}$ denote the closure in $X$ of the
image of $g\hH^2$, which identifies with the Hilbert modular surface
attached to $F_j$ and the congruence subgroup $g^{-1}\Gamma g \cap
{\rm SL}(2, \Oo_{F_j})$. It is the natural analogue of the
Hirzebruch-Zagier cycle on a Hilbert modular surface.

\medskip

The proof of Theorem A has as a consequence the following:

\medskip

\noindent{\bf Theorem B} \, \it Let $K$ be biquadratic with
intermediate quadratic fields $F_1, F_2, F_3$. Define $N>0$ by
$N\Z=\nN\cap\Z$, and assume that the modular curve $X_0(N)$ has
genus $>0$. Then there exist $g_1,g_2 \in {\rm GL}_2^+(K)$ such that
$Y_{F_1,g_1}, Y_{F_2,g_2}$ span a $2$-dimensional subspace of
$r_{{\rm alg},\Q}$. Consequently,
$$
{\rm dim}_\Q CH^2(X) \, \geq 2.
$$
\rm

\medskip

By the product structure, one gets non-trivial, decomposable classes
in $CH^3(X,1)$. A refinement will be discussed in section 11.

\bigskip

\section{Contribution from the boundary}

\medskip

By the {\it decomposition theorem}, there is a short exact sequence
$$
0\to IH^4(Y^\ast) \to H^4(X) \to H^4_{\tilde Y^\infty}(X) \to 0,
$$
both as Galois modules and as $\Q$-Hodge structures. Here $IH^\ast$
is the Goresky-MacPherson's middle intersection cohomology. We
obtain, for $\alpha = \{\ell,k\}$, an, alg,
$$
r_\alpha(X) \, = \, r_\alpha(Y^\ast)+ r_\alpha^\infty,
$$
where $r_\alpha(Y^\ast)$, resp. $r_\alpha^\infty$, is the
$\alpha$-rank associated to $IH^4(Y^\ast)$, resp. $H^4_{\tilde
Y^\infty}(X)$.

\medskip

The cohomology with supports in $\tilde Y^\infty$ has a nice
description:

$$
H^\ast_{\tilde Y^\infty}(X)  \, \simeq  \, \oplus_{\sigma: {\rm
cusp}} \, \, H^\ast_{D_\sigma}(X),
$$
where $D_\sigma$ is a divisor with normal crossings (DNC) with
smooth irreducible components $D^i_\sigma$. If $D^{i,j}_\sigma$
denotes $D^i_\sigma \cap D^j_\sigma$, there is an exact sequence
$$
\sum\limits_i H^2(D^i_\sigma)(1)  \to H^4_{D_\sigma}(X)(2) \to
\sum\limits_{i \ne j^{}} H^0(D^{i,j}_\sigma)
$$
Since $D^i_\sigma$ is toric, its $H^2$ is generated by divisors.
This implies the following string of equalities for large $k$:
$$
r_{\rm alg,k}^\infty = r_{\ell,k}^\infty = r_{{\rm an}, k}^\infty =
r_{\rm Hg}^\infty.
$$
All but the last equality on the right remain in force for any
number field $k$.

\medskip

Hence the problem reduces to understanding the $r_\alpha(Y^\ast)$
for various $\alpha$ and explicating their relationships with each
other.

\bigskip

\section{The action of Hecke correspondences}

\medskip

If $g \in $GL$_2^+(K)$, there are two maps $Y_{\Gamma_g}\to
Y_\Gamma$, with $\Gamma_g = \Gamma \cap g^{-1}\Gamma g$, inducing an
algebraic correspondence $T_g$, which extends to $Y_\Gamma^\ast$.
The algebra ${\mathcal H}$ of such Hecke correspondences acts
semisimply on cohomology. This leads to a ${\mathcal H} \times
{\mathcal G}_\Q$-equivariant decomposition
$$
IH^4(Y^\ast) \, \simeq \, V_{\rm res} \oplus V_{\rm cusp}
$$
where the submotive $V_{\rm res}$ is algebraic, and the cuspidal
submotive $V_{\rm cusp}$ is the ``interesting'' part (see below).

To be precise, the {\it residual part} $V_{\rm res}$ is spanned by
the intersections of Chern classes of certain universal line bundles
${\mathcal L}_{ij}, 1\leq i,j\leq 4$ occurring at every level $\nN$.
In the complex realization they are defined by the SL$(2,
\R)^4$-invariant differential forms $\eta_i \wedge\eta_j$, $1 \leq i
\ne j \leq 4$ on $\hH^4$, where for each $j$, $\eta_j={\rm
pr}_j^\ast(dz_j \wedge d{\overline z}_j)$, with pr$_j: \hH^4 \to
\hH$ being the $j$-th projection.

\medskip

\noindent{\bf Remark}: \, In our case one knows enough about the
Galois modules which occur in $H^\cdot(X)$ to be able to get a
direct sum decomposition:
$$
H^\cdot(X) \, \simeq \, IH^\cdot(Y^\ast) \oplus  H^\cdot_{\tilde
Y^\infty}(X)
$$
$IH^\cdot(Y^\ast)$ can be cut out as a direct summand by an
algebraic cycle modulo homological equivalence. (This has recently
been done for general Shimura varieties by A.~Nair ([N]).) The
reason is that the Hecke correspondences act on an inverse limit of
a family of toroidal compactifications of $Y$, though not on any
individual one. However, $IH^\cdot(Y^\ast)$ is not immediately a
Chow motive, since the algebra of Hecke correspondences modulo
rational equivalence is not semisimple.

\medskip

The K\"unneth components of $\Delta$ in $IH^8(Y^\ast\times Y^\ast)$
are algebraic, since it is known that $IH^1(Y^\ast)=IH^3(Y^\ast)=0$
and $IH^2(Y^\ast)$ is algebraic, being purely of Hodge type $(1,1)$.

\medskip

There is a further ${\mathcal H} \times {\mathcal G}_\Q$-equivariant
decomposition:
$$
V_{\rm cusp} \, = \, \oplus_\varphi V(\varphi)^{m(\varphi)},
$$
where $\varphi$ runs over holomorphic Hilbert modular cusp forms of
level $\nN$, which have (diagonal) weight $2$, and $m(\varphi)$ is a
certain multiplicity which is $1$ if $\varphi$ is a {\it newform},
i.e., not a cusp form of level a proper divisor of $\nN$.

\medskip

\section{ The submotives of rank ${\bf 16}$}

\medskip

It is now necessary to understand $V(\varphi)$ for a Hilbert modular
newform $\varphi$ of weight $2$ and level $\nN$. It is easy to see
that $V_B(\varphi)$ is $16$-dimensional, generated over $\C$ by the
differential forms $\varphi(z)dz_I\wedge d{\overline z}_J$ of degree
$4$ for partitions $\{1,2,3,4\}=I\cup J$.

Let $\pi$ be the cuspidal automorphic representation of GL$(2,
\A_K)$ of trivial central character associated to $\varphi$. We will
write $V(\pi)$ instead of $V(\varphi)$. By R.L.~Taylor and
Blasius-Rogawski, one can associate a $2$-dimensional irreducible
representation $W_\ell(\pi)$ of ${\mathcal G}_K$ such that
$L(s,W_\ell(\pi))=L(s,\pi)$, i.e., $\forall \mathcal P \nmid
\mathcal N$,
$$
{\rm tr}\left(Fr_\mathcal P\vert W_\ell(\pi)\right) \, = \,
a_\mathcal P(\pi).
$$

By Brylinski-Labesse, as refined by Blasius,
$$
{\rm tr}\left(Fr_\mathcal P\vert V_\ell(\pi)\right) \, = \, L(s,
\pi_\mathcal P; r), \, \, \, \forall \mathcal P\nmid \mathcal N.
$$
Since $V_\ell(\pi)$ and $W_\ell(\pi)$ are semisimple, we get the
following isomorphism by Tchebotarev:
$$
V_\ell(\pi)\vert_{{\mathcal G}_K}\simeq\otimes_{\tau\in{\rm
Gal}(K/\Q)}W_\ell(\pi)^{[\tau]}
$$
To be precise, $V_\ell(\pi)$ is the {\it tensor induction} ([CuR])
of $W_\ell(\pi)$ from $K$ to $\Q$:
$$
V_\ell(\pi) \simeq {}^{\otimes}{\rm Ind}_{{\mathcal G}_K}^{{\mathcal
G}_\Q}(W_\ell(\pi))
$$
This identity makes it possible to compute the Tate classes.

We need the following result:

\medskip

\noindent{\bf Theorem C} ([Ra3,4]) \, \it The $L$-function of
$V_\ell(\pi)$ admits a meromorphic continuation to the whole
$s$-plane with a functional equation of the form
$$
L(s, V_\ell(\pi)^\vee) = \varepsilon(s,
V_\ell(\pi))L(5-s,V_\ell(\pi)),
$$
where $\varepsilon(s, V_\ell(\pi))$ is an invertible function on
$\C$ and the superscript $\vee$ indicates the dual. \rm

\medskip

From this one also gets the analogous statement about the
$L$-function of $X$, which is a product of these $L(s,V_\ell(\pi))$
with an abelian $L$-function.

\bigskip

\section{Strategy for algebraicity}

\medskip

When $r_\ell(\pi) \ne 0$, we show first that $L(s,V_\ell(\pi))$ has
a pole at $s=3$. But then we also show, using a specific integral
representation, that for a suitable quadratic field $F\subset K$,
the function $L_{1,F}(s):= L(s, {}^{\otimes}{\rm Ind}_{{\mathcal
G}_K}^{{\mathcal G}_F}(W_\ell(\pi))\otimes \nu)$ has a {\it simple}
pole. What we do then is to construct an algebraic cycle ${}^\nu
Z_F$, and prove $\int_{{}^\nu Z_F} \omega \ne 0$ for a $(2,2)$-form
$\omega$ by realizing it as res$_{s=2}L_{1,F}(s)$. Using the
previous section, we first prove

\medskip

{} \noindent{\bf Proposition} \, \it For any Dirichlet character
$\chi$,
$$
r_\ell(\pi, \chi)=r_{\rm an}(\pi, \chi)\leq 2
$$
$r_\ell(\pi,\chi)=1$ iff a twist $\pi\otimes\nu$ is fixed by an
involution $\tau \in $Gal$(K/\Q)$, while $r_\ell(\pi,\chi)=2$ iff
$K$ is biquadratic and $\pi\otimes\nu$ is Gal$(K/\Q)$-invariant.\rm

\medskip

Now the problem is to construct algebraic cycles in the
$(\pi\otimes\nu)$-eigenspace, which are of infinite order modulo
homological equivalence, to account for these poles when
$\pi\otimes\nu$ is fixed by one or more (non-trivial) elements of
Gal$(K/\Q)$.

\bigskip

\section{Cycles}

\medskip

Let $F$ be a quadratic subfield of $K$, with corresponding embedding
${\mathfrak h}^2\hookrightarrow{\mathfrak h}^4$. Recall that if
$g\in{\rm GL}_2^+(K)$, the image of the translate $g{\mathfrak h}^2$
under the projection ${\mathfrak h}^4\to\Gamma\backslash{\mathfrak
h}^4$ defines a surface $\Delta_g\backslash {\mathfrak h}^2$ with
closure $Z_{F,g}$ in $X$. This is an example of a Hirzebruch-Zagier
cycle; see [K] for a definition of such cycles for orthogonal
Shimura varieties.

For any abelian character $\mu$, there is a $\mu$-{\it twisted
Hirzebruch-Zagier cycle} $Z_{F,g}^\mu$ of codimension $2$ in $X$.
This is defined (cf. [Ra4]) by composing the above construction with
a {\it twisting correspondence} defined by $\mu$, which sends, for
every $\pi$, the $\pi$-isotypic subspace onto the $(\pi\otimes
\mu)$-isotypic subspace of the cohomology.

\medskip

Suppose $r_\ell(\pi,\chi) (=r_{\rm an}(\pi,\chi))$ is $>0$. Then
there is a quadratic subfield $F$ and a cusp form $\pi_1$ on
GL$(2)/F$ such that $\pi\otimes\nu \simeq(\pi_1)_K$. As one would
hope, a twisted Hirzebruch-Zagier cycle $Z_{F,g}^\mu$ provides the
requisite algebraic cycle to get the Tate conjecture.

There is a real {\bf subtle point} here which separates it from the
work of Harder, Langlands and Rapoport ([HLR]) on the divisors on
Hilbert modular surfaces:

{\it The period of a $(2,2)$-form on $X$ (defined by $\pi$) over
$Z_{F,g}^\mu$ is non-zero, but it is the residue of a {\bf different
$L$-function}, namely $L_{1,F}(s)$, which {\bf does not divide}
$L(s,V_\ell(\pi))$}!

The residues of the two $L$-functions are presumably related in a
non-trivial way, but this is not known. It is an intriguing problem
to try to understand this better.

\bigskip

\section{Hodge classes}

\medskip

By hypothesis, we need only consider those $\pi$ which are of level
$\nN \ne \Oo_K$. We may then fix a prime divisor ${\mathcal P}$ of
$\nN$ and consider a quaternion algebra $B/K$ which is ramified only
at three infinite places and at ${\mathcal P}$. By the
Eichler-Shimizu-Jacquet-Langlands correspondence, there exists a
corresponding cusp form $\pi'$ on $B^\ast$ giving rise to a
submotive $V(\pi')$ of $H^4(R_{K/\Q}(C))$ for a Shimura curve
$C=\Delta\backslash {\mathfrak h}$ defined over $K$. As
Gal$(\overline \Q/K)$-modules,
$$
V_\ell(\pi') \, \simeq \, V_\ell(\pi),
$$
implying that
$$
r_{\ell,E}(\pi') \, = \, r_{\ell, E}(\pi),\leqno(i)
$$
for any number field $E \supset K$.

On the Hodge side we need the following result, proved jointly with
V.K.~Murty, which is really a statement about periods:

\medskip

\noindent{\bf Theorem} ([Mu-Ra2]) \, \it As $\Q$-Hodge structures,
$$
V_B(\pi') \, \simeq \, V_B(\pi).
$$
\rm

\medskip

The proof compares the coefficients of the Shimura liftings of $\pi$
and $\pi'$ to forms of weight $3/2$, which live on the two-sheeted
covering group of GL$(2)/K$.

\medskip

By Deligne ([DMOS]), every Hodge class on an abelian variety gives
rise to a Tate class. One gets from this the equality over a
sufficiently large field $E\supset K$:
$$
r_{\ell,E}(\pi') \, = \, r_{\rm Hg}(\pi').\leqno(ii)
$$

On the other hand, by the Theorem with Murty, one also gets
$$
r_{\rm Hg}(\pi')=r_{\rm Hg}(\pi).\leqno(iii)
$$
Combining (i), (ii) and (iii), we see that the Hodge conjecture for
$V_B(\pi)$ follows from the Tate conjecture for $V_\ell(\pi)$ over
$E$.

\bigskip

\section{Where the cycles come from}

\medskip

The method of proof furnishes the following:

\medskip

\noindent{\bf Theorem B$^\prime$} \, \it Let $F$ be a quadratic
subfield of $K$. Then a twisted Hirzebruch-Zagier cycle of
codimension $2$ on $X$ associated to $F$ contributes to $V(\pi)$ iff
a twist of $\pi$ is a base change from $F$. When $K$ is biquadratic
and a twist of $\pi$ is base changed from $\Q$, i.e., attached to an
elliptic cusp form $h$, we get
$$
r_{{\rm alg},\Q}(\pi) \, = \, 2,
$$
with the Hecke twisted Hilbert modular surfaces from two subfields
$F_1, F_2$, say, give non-trivial independent algebraic classes of
codimension $2$. \rm

\medskip

As a consequence, the decomposable part of $CH^3(X_\pi,1)$ has rank
at least $2$ when $\pi$ is a base change from $\Q$ and $K$
biquadratic. (Here $X_\pi$ refers to the submotive of $[X]$ in
degree $4$ cut out by $\pi$.) But since $\Z^\ast = \{\pm 1\}$, such
classes will not come from a regular, proper model ${\mathcal X}$ of
$X$ (when such a model exists).

\medskip

Now let $N>0$ be defined by $N\Z=\nN\cap \Z$. Then when the modular
curve $X_0(N)$ has positive genus, there exists at least one
(elliptic) newform $h$ of weight $2$ and level $N$, and the base
change to $K$ of the associated cuspidal automorphic representation
on GL$(2)/\Q$, assures the existence of a $\pi$ as in Theorem
B$^\prime$. Hence Theorem B follows from Theorem B$^\prime$.

\bigskip

\section{What to expect}

\medskip

A straight-forward calculation shows that
$$
-{\rm ord}_{s=2}L_\infty(s, V(\pi)) \, = \, 3.
$$

In the {\it biquadratic} case, if a twist of $\pi$ is a base change
from $\Q$, we have
$$
r_{\rm an,\Q}(\pi^\vee) \, = \, 2.
$$
Hence by the functional equation,
$$
{\rm ord}_{s=2} L(s, V(\pi)) \, = \, 1.
$$
In this case, Beilinson predicts (in [Be]) the existence of a
non-trivial class $\beta$ in $CH^3(X_\pi,1)_\Q$ which comes from a
proper model ${\mathcal X}$ over $\Z$, to account for this simple
zero of the $L$-function.

\medskip

In the {\it cyclic} case, $r_{\rm an,\Q}(\pi^\vee)=1$, and so the we
should have two independent classes in the higher Chow group. The
general philosophy is that it is much harder to produce classes in a
motivic cohomology group (or a Selmer group) which is conjecturally
of rank bigger than one, and this is why we are not at present
concentrating on this (cyclic) case.

\bigskip

\section{ Elements in ${\bf CH^3(X,1)_\Q}$}

\medskip

Let $K/\Q$ be biquadratic with quadratic subfields $F_1, F_2, F_3$.
For $g_1, g_2, g_3 \in {\rm GL}_2^+(K)$, consider the surfaces
$Z_i=Z_{F_i,g_i}, 1\leq i\leq 3$, which are Hecke translates of the
three Hilbert modular surfaces in $X$ associated to $\{F_i\vert1\leq
i\leq 3\}$. Put
$$
C_{i,j} \, = \, Z_i\cap Z_j, \quad {\rm for} \quad 1\leq i\ne j\leq
3.
$$

\medskip

\noindent{\bf Theorem C} \, \it One can find $g_1,g_2,g_3$ in ${\rm
GL}_2^+(K)$, and an integer $m>0$, such that, up to modifying the
construction by decomposable classes supported on the boundary and
on $X_{\rm res}$, we have for each $i \leq 3$ and a permutation
$(i,j,k)$ of $(1,2,3)$,
$$
m(C_{i,j}-C_{i,k}) \, = \, {\rm div}(f_i)
$$
for some functions $f_i$ on $Z_i$. \rm

\medskip

Here is the {\bf basic idea}: Each $Z_i$ is a Hecke translate of a
Hilbert modular surface $S_i$, which is simply connected. So
$Pic(Z_i)= NS(Z_i)$. We show that the homology class
$[g_i^{-1}(C_{i,j}-C_{i,k})]$ is trivial in $H^2(S_i,\Q(1))$, up to
modifying by the {\it trivial classes} coming from the boundary and
the residual part. Thanks to the explicit description of the
algebraic cycles (modulo homological equivalence) on Hilbert modular
surfaces (see [HLR], [MuRa1]), one knows that in the situation we
are in, the divisor classes are spanned by Hecke translates of
modular curves. From this it is not too difficult to show (for
suitable $\{g_i\}$) the homological triviality (modulo trivial
cycles) of $[g_i^{-1}(C_{i,j}-C_{i,k})]$ {\it when} the Hilbert
modular surface $S_i$ has geometric genus $1$; in fact it is enough
to know that there is a unique base changed newform $\pi_i$ of
weight $2$ over $F_i$ for this. In general one has to deal with
several newforms, and one uses a delicate refinement of an argument
of Zagier ([Z], page 243). The subtlety comes from the fact that one
needs to deal with three quadratic fields at the same time.

\medskip

A simple example is when $g_1^{-1}g_2, g_1^{-1}g_3$ are diagonal
matrices in GL$_2^+(K)$, not in GL$_2^+(\Q)$, fixed by Gal$(K/F_1)$,
Gal$(K/F_3)$ respectively. For each $i\leq 3$, $C_{ij}, C_{ik}$ are
Hecke translates of modular curves on $Z_i$.

\medskip

Thanks to Theorem C, the formal sum
$$
\sum\limits_{i=1}^3 \, (Z_i,f_i)
$$
satisfies $\sum_i {\rm div}(f_i) = 0$ as a codimension $2$ cycle on
$X$, and hence defines a class in $CH^3(X,1)\otimes \Q$.

\medskip

\noindent{\bf Problem}: \, \it Compute, for $\omega\in
H^{(2,2)}(X(C))$,
$$
\sum_i \int_{Z_i} \log\vert f_i\vert \, \omega
$$
\rm

\medskip

We can understand this period integral a bit better in the analogous
situation where $X$ is the four-fold product of a modular curve
$X_0(N)$ for prime level $N$, the simplification arising from the
fact that one can reduce to considering $f$ of the form
$\sum_{j=1}^r f_{1,j}\otimes f_{2,j}\otimes f_{3,j}\otimes f_{4,j}$,
with each $f_{i,j}$ a translate of a modular unit on $X_0(N)$. Here
the Hilbert modular surfaces are replaced by the twisted images of
$X_0(N)^2 \to X_0(N)^4$, which are not simply connected, but anyhow,
their Pic$^0$ is generated by elementary divisors of degree $0$, and
this suffices. Work is in progress to prove for $N=11$, by a
combination of theoretical and numerical arguments, that the
integral is non-zero for a suitable choice; this will show that the
class is of infinite order. We hope to investigate if such a class
can give any information on the Bloch-Kato Selmer group at $s=0$ of
the Sym$^4$ motive, twisted by $\Q(2)$, of $X_0(N)$.

\bigskip

\section{The modular complex}

\medskip

Let $X$ be a smooth, toroidal compactification of a Shimura variety
of dimension $n$ over its natural field of definition $k$. Consider
the class of closed irreducible subvarieties of $X$, called {\it
modular}, generated by components of the twisted Hecke translates of
Shimura subvarieties, the boundary components, the components of
their intersections, and so on. For every $p \geq 0$, let $X_{\rm
mod}^p$ denote the set of such subvarieties of codimension $p$. The
{\bf modular points} in this setting will be the CM points and the
points arising from successive intersections of components at the
boundary. To give a concrete example, consider the case of a Hilbert
modular surface $X$ which is obtained by blowing up each cusp into a
cycle of rational curves. When $\Gamma_0(\nN)$ is torsion-free, the
{\it modular points} on $X$ will be the CM points and the points
where the rational curves over cusps intersect.

We can now consider the {\bf modular analogue} of the Gersten
complex, namely
$$
\dots \to \coprod\limits_{W\in X_{\rm mod}^{p+2}} K_2(k(W)) \to
\coprod\limits_{Z\in X_{\rm mod}^{p+1}} k(Z)^\ast \to
\coprod\limits_{Y\in X_{\rm mod}^p} \Z.Y
$$
We may look at the homology of this complex, and denote the
resulting groups - the first two from the right -- by $B^p(X)$ and
$B^{p+1}(X,1)$. There are natural maps
$$
B^p(X) \to CH^p(X), \quad {\rm and} \quad B^{p+1}(X,1) \to
CH^{p+1}(X,1)
$$
Denote the respective images by $CH^p_{\rm mod}(X)$ and
$CH^{p+1}_{\rm mod}(X,1)$. The nice thing about these groups is that
since the modular subvarieties are all defined over number fields,
the building blocks do not change from $\overline \Q$ to $\C$. It
will be very interesting (exciting?) to try to verify, in some
concrete cases of dimension $\geq 2$, whether $CH^p_{\rm mod}(X)$
and $CH^{p+1}_{\rm mod}(X,1)$ are finitely generated over $k$.

Note that the classes we consider in this article are {\it modular}
in this sense. When $n=2m$ and when $L(s, M)$ has, for a simple
submotive $M$ of $H^n(X)$, a {\it simple zero} at $s=m$, resp. a
{\it simple pole} at $s=m+1$, it is tempting to ask, in view of the
known examples, if there is an element of $CH^{m+1}_{\rm mod}(X,1)$,
resp. $CH^m_{\rm mod}(X)$, which {\it explains} it. When $n=2m-1$
and $L(s, M)$ has a simple zero at $s=m$, one could again ask if it
is explained by a {\it modular} element of $CH^m(X)^0$, the
homologically trivial part. For modular curves $X_0(N)/\Q$, one has
striking evidence for this in the work of Gross and Zagier. For
Hilbert modular surfaces it is again true ([HZ], [HLR]). The
situation is {\it not} the same when the order of pole or zero is
$2$ or more, especially over non-abelian extensions of $k$
([MuRa1]).

\bigskip

\section*{\bf Bibliography}

\begin{description}

\item[{[Be]}] A.A.~Beilinson, Higher regulators and values of
$L$-functions, Journal of Soviet Math. {\bf 30}, No. 2, 2036--2070
(1985).

\item[{[B]}] S.~Bloch,  Algebraic cycles and higher $K$-theory.
Advances in Math. {\bf 61}, no. 3, 267--304 (1986).

\item[{[CuR]}] C.W.~Curtis and I.~Reiner, Methods of
representation theory I, Wiley, NY (1981).

\item[{[DMOS]}] P.~Deligne, J.S.~Milne, A.~Ogus and K-Y.~Shih,
Hodge cycles, motives, and Shimura varieties. Lecture Notes in
Mathematics {\bf 900}. Springer-Verlag, Berlin-New York (1982).

\item[{[HLR]}] G.~Harder, R.P.~Langlands and
M.~Rapoport, Algebraische Zykeln auf Hilbert-Blumenthal-Fl\"achen,
Crelles Journal {\bf 366} (1986), 53--120.

\item[{[HZ]}] F.~Hirzebruch and D.~Zagier, Intersection numbers of
curves on Hilbert modular surfaces and modular forms of Nebentypus,
Inventiones Math. {\bf 36}, 57--113  (1976).

\item[{[K]}] S.~Kudla, Algebraic cycles on Shimura varieties of
orthogonal type, Duke Math. Journal  {\bf 86}, no. 1, 39--78 (1997).

\item[{[MuRa1]}] V.K.~Murty and D.~Ramakrishnan, Period relations and
the Tate conjecture for Hilbert modular surfaces, Inventiones Math.
{\bf 89}, no. 2 (1987), 319--345.

\item[{[MuRa2]}] V.K.~Murty and D.~Ramakrishnan, Comparison of $\Q$-Hodge
structures of Hilbert modular varieties and Quaternionic Shimura
varieties, {\it in preparation}.

\item[{[N]}] A.~Nair, Intersection cohomology, Shimura varieties,
and motives, preprint (2003).

\item[{[Ra1]}] D.~Ramakrishnan, Arithmetic of Hilbert-Blumenthal
surfaces, CMS Conference
Proceedings {\bf 7}, 285--370 (1987).

\item[{[Ra2]}] D.~Ramakrishnan, Periods of integrals arising from
$K_1$ of Hilbert-Blumenthal surfaces, preprint (1988); and  Valeurs
de fonctions $L$ des surfaces d'Hilbert-Blumenthal en $s=1$, C. R.
Acad. Sci. Paris S\'er. I Math. {\bf 301}, no. 18, 809--812 (1985)

\item[{[Ra3]}] D.~Ramakrishnan, Modularity of solvable Artin
representations of GO$(4)$-type, International Mathematics Research
Notices (IMRN) {\bf 2002}, No. 1 (2002), 1--54.

\item[{[Ra4]}] D.~Ramakrishnan, Algebraic cycles on Hilbert modular
fourfolds and poles of $L$-functions, in {\it Algebraic Groups and
Arithmetic},  221--274, Tata Institute of Fundamental Research,
Mumbai (2004).

\item[{[Sch]}] A.J.~Scholl, Integral elements in $K$-theory and products
of modular curves, in {\it The arithmetic and geometry of algebraic
cycles}, 467--489, NATO Sci. Ser. C Math. Phys. Sci. {\bf 548},
Kluwer Acad. Publ., Dordrecht {\bf 2000}.

\item[{[Z]}] D.~Zagier, Modular points, modular curves, modular
surfaces and modular forms, {\it Workshop Bonn 1984}, 225--248,
Lecture Notes in Math. {\bf 1111}, Springer, Berlin-New York (1985).

\bigskip

\end{description}

\vskip 0.3in

\noindent Dinakar Ramakrishnan

\noindent Department of Mathematics

\noindent California Institute of Technology, Pasadena, CA 91125.

\noindent dinakar@its.caltech.edu

\end{document}